\documentclass[letterpaper,11pt]{article}
\usepackage[letterpaper, left=1in,right=1in,top=1in,bottom=1in]{geometry}   % Use for the final submission
\usepackage{color}
\usepackage{amssymb}
\usepackage{epsfig}
\usepackage{amsmath}
\usepackage{amsfonts}
\usepackage{graphicx}
\usepackage{subfigure}
\usepackage{color}
\usepackage{graphicx}
\usepackage{epsfig}
\usepackage{cite}
\usepackage{amsmath}
\usepackage{amssymb}
\usepackage{mathrsfs}
\usepackage{amsfonts}
\usepackage{dsfont}
\usepackage{graphicx}
\usepackage{mathrsfs}
\usepackage{hyperref}
\usepackage{hyperref}
\usepackage{titlesec}
\usepackage[T1]{fontenc}
\usepackage{ragged2e}
\justifying

\DeclareGraphicsExtensions{.pdf,.png,.jpg,.JPG}
\usepackage{lastpage}
\usepackage{fancyhdr}
\raggedright
\titleformat{\section}{
	\vspace{12pt}\scshape\raggedright\large
}{}{0em}{\textsc}

\begin{document}
	\thispagestyle{empty}
	
	\begin{tabular*}{\textwidth}{l@{\extracolsep{\fill}}r}
		\textbf{\centerline{Combining Learning and Control for Data-driven Approaches of Cyber-Physical Systems}}\vspace{5pt}\\ 
		\centerline{Andreas A. Malikopoulos}\\
        \centerline{Cornell University}\\
        \centerline{amaliko@cornell.edu}\\
	\end{tabular*}

\vspace{+5EX}

\centerline{\textbf{Abstract}}

\justify{

Cyber-physical systems (CPS), in most instances, represent systems of systems with an informationally decentralized structure such as emerging mobility systems, networked control systems, sustainable manufacturing,  smart power grids, power systems, mobility markets,  social media platforms, cooperation of robots, and internet of things. 
To optimize the operation of such systems, we typically assume an ideal model. Such model-based control approaches cannot effectively facilitate optimal solutions with performance guarantees due to the discrepancy between the model and the actual CPS. On the other hand, in most CPS there is a large volume of data with a dynamic nature which is added to the system gradually in real time and not altogether in advance. Thus, traditional supervised learning approaches cannot always facilitate robust solutions using data derived offline. By contrast, applying reinforcement learning approaches directly to the actual CPS might impose significant implications on the safety and robust operation of the system. 
The \textbf{overarching goal} of the Information and Decision Science (IDS) Lab is to investigate how to circumvent these challenges by developing \textbf{data-driven} system approaches at the intersection of \textbf{learning} and \textbf{control}.	The emphasis is on how to improve \textbf{energy efficiency} and reduce \textbf{greenhouse gas emissions} in applications related to emerging mobility systems, e.g., connected and automated vehicles (CAVs), shared mobility, sociotechnical systems, and smart cities, and thus contribute to the \textbf{health of the planet}.

\section{\textbf{Self-Learning Powertrain Control}}
\justify{
	
My interest in developing control algorithms that could make systems able to learn their optimal operation started early on while I was still at graduate school when I read an article about the discrepancy between the true fuel economy of a vehicle and the one posted on the window sticker. The article discussed the implications of the driver’s driving style on engine operation and stated that the state-of-the-art control methods, at that time, consisted of static controllers that could not optimize engine operation for different driving styles but only for predetermined ones. This article provided inspiration that eventually led to the formation of my dissertation topic. In my dissertation \cite{Malikopoulos2011b}, I developed the theoretical framework \cite{Malikopoulos2007,Malikopoulos2010a,Malikopoulos2009e,Malikopoulos2009} and control algorithms  \cite{Malikopoulos2007a,Malikopoulos2009a,Malikopoulos2006} that can turn the engine of a vehicle into an autonomous intelligent system capable of learning its optimal operation in real time while the driver is driving the vehicle. I modeled the evolution of the state of the engine as a control Markov chain \cite{Malikopoulos2010b} and proved \cite{Malikopoulos2009b} that it eventually converges to a stationary probability distribution deemed characteristic of the driver’s driving style. Through this approach, the engine progressively perceives the driver’s driving style \cite{Malikopoulos2008b} and eventually learns to operate in a manner that optimizes specified performance criteria, e.g., fuel economy and emissions with respect to the driver’s driving style. The framework also allows the engine to identify the driver, and thus it can adjust its operation to be optimal for any driver based on what it has learned in the past regarding her/his driving style. The outcome of my dissertation research eventually led to a US patent \cite{Malikopoulos2013d}.

Moving to General Motors Research \& Development as a Senior Researcher, I had the chance to continue working on self-learning control for advanced powertrain systems. I led several projects on autonomous intelligent propulsion systems and developed computational mathematical models and control algorithms for making highly energy-efficient and eco-friendly vehicles. I was a member of the team that successfully demonstrated the implementation of self-learning control algorithms \cite{Malikopoulos2014c} in two demo vehicles, Saturn Aura and Opel Vectra. 
}

\section{\textbf{Hybrid-Electric Vehicles (HEVs) and Plug-in HEVs}}
\justify{
When I joined Oak Ridge National Laboratory (ORNL), I had the chance to work across different technical areas, including stochastic optimal control \cite{Malikopoulos2011a,Malikopoulos2015,Malikopoulos2015b}, optimal design and power management control and routing of hybrid electric vehicles (HEVs) and plug-in HEVs (PHEVs) 
\cite{Park2011,Malikopoulos2011,Malikopoulos2013b,Malikopoulos2013,Malikopoulos2014,Malikopoulos2014a,Shaltout2014,Pourazarm2014, Shaltout2015b,Malikopoulos2015d,Malikopoulos2016b}, and driver’s feedback systems \cite{Rios-Torres2016a,Malikopoulos2012,Malikopoulos2013a}. The latter eventually led to a technology \cite{Malikopoulos2015a} that was licensed in SanTed Project Management LLC. I also contributed to the solution of problems that included smart buildings aimed at optimizing energy system parameters to (1) improve sustainability, (2) facilitate cost-effective energy generation, and (3) allocate demand optimally to different energy sources, e.g., solar, renewable, etc \cite{Sharma2016,Dong2017a,Dong2016}. On the fundamental research front, I established the theoretical framework for the analysis and stochastic control of complex systems consisting of interactive subsystems \cite{Malikopoulos2016}. In particular, I developed a duality framework and showed that the Pareto control policy minimizes the long-run expected average cost criterion of the system while also presented a geometric interpretation of the solution and conditions for its existence. I provided theoretical results showing that the Pareto control policy provides an equilibrium operating point among the subsystems, and if the system operates at this equilibrium, then the long-run expected average cost per unit time is minimized. This result implies that the Pareto control policy can be of value when we seek to derive the optimal control policy for complex systems online. Later on, and in my role as the Deputy Director of the Urban Dynamics Institute at ORNL, I developed several initiatives with the goal to investigate how we can use scalable data and informatics to enhance understanding of the environmental implications of CAVs and improve transportation sustainability and accessibility. I contributed towards the development of a decentralized optimal control framework whose closed-form solution exists under certain conditions, and which, based on Hamiltonian analysis, yields for each vehicle the optimal acceleration/deceleration, in terms of fuel consumption. The solution allows the vehicles to cross merging roadways without creating congestion, and under the hard safety constraint of collision avoidance \cite{Rios-Torres2015,Rios-Torres2017b,Rios-Torres2017a,Rios-Torres2017,Rios-Torres2016b}.

}

\section{\textbf{Connected and Automated Vehicles}}
\justify{

Emerging mobility systems are typical cyber-physical (CPS) systems where the cyber component (e.g., data and shared information through vehicle-to-vehicle and vehicle-to-infrastructure communication) can aim at optimally controlling the physical entities (e.g., CAVs, non-CAVs). A mobility system encompasses the interactions of three heterogeneous dimensions: (1) transportation systems and modes, e.g., CAVs, shared mobility, and public transit integrated with advanced control algorithms, (2) social behavior of drivers, operators (for autonomous vehicles), and travelers (or pedestrians) interacting with these systems, and (3) information management of data available and shared information. The constellation of these three dimensions constitutes a sociotechnical system that should be analyzed holistically. The CPS nature of emerging mobility systems is associated not only with technological and information management dimensions but also with human adoption (social dimension). My students and I, in conjunction with my collaborators, have made contributions on the technological dimension of mobility systems by developing control algorithms for optimal coordination of CAVs \cite{Malikopoulos2018,Malikopoulos2018b,Rios-Torres2018,Stager,Zhao2018,Zhang2017a,Zhang,Zhao2018a, assanis2018ITSC,Beaver2020DemonstrationCity,chalaki2019optimal,chalaki2020experimental,chalaki2020hysteretic,chalaki2020ICCA,chalaki2020TCST,chalaki2020TITS, Connor2020ImpactConnectivity,jang2019simulation,Mahbub2019ACC,Mahbub2020ACC-1,mahbub2020ACC-2,mahbub2020Automatica-2,mahbub2020decentralized,mahbub2020sae-1,mahbub2020sae-2,malikopoulos2019ACC,Malikopoulos2019CDC,Malikopoulos2020,Zhao2018CTA,Zhao2018ITSC,Zhao2019CCTA-1,zhao2019CCTA-2,Ray2021DigitalCity,chalaki2020hysteretic,Sumanth2021,chalaki2021Reseq,K122021AExperiment,chalaki2021CSM,mahbub2021_platoonMixed,mahbub2020Automatica-2,mahbub2022ACC,chalaki2021CSM,K122021AExperiment,Ray2021DigitalCity} and identifying potential research paths with connected autonomous systems \cite{zhao2019enhanced}. 
The framework has been extended to expand performance and the feasibility domain for CAVs at signal-free intersections\cite{tzortzoglou2023approach,tzortzoglou2023performance,tzortzoglou2024feasibility}.

Given that a transportation network with a 100\% penetration rate of connected and automated vehicles (CAVs) is not expected to be realized in the next few decades, addressing planning, control, and navigation for CAVs in mixed traffic, considering the co-existence of human-driven vehicles (HDVs) with various driving styles, is imperative.
Our early work on that topic focused on mixed-traffic platoon formation \cite{mahbub2022_ifac,mahbub2023_automatica} which aims to guarantee string stability between CAVs and HDVs. 
However, we are currently more interested in problems related to interaction-driven traffic scenarios, such as merging at roadways and roundabouts, crossing intersections, and lane changing.
In these scenarios, instead of stabilizing the operation of the CAVs, the controllers must optimize the CAVs' behavior to ensure safety, improve travel time, avoid gridlocks, and minimize traffic disruption and human discomfort. 
One typical control technique that can be used is model predictive control (MPC), due to its flexibility in considering multiple objectives and constraints, as well as its robustness to uncertainty caused by human drivers due to the receding horizon implementation.
While there is an increasing number of research articles on MPC for individual CAVs in mixed traffic, a gap in the existing research is the ability to adapt the control designs given different HDVs' driving styles.
Generally, non-adaptive MPC with fixed weights cannot guarantee to work well in this application.
We developed a control framework \cite{Le2022CDC} to address the motion planning problem for single CAVs while interacting with HDVs, which combines three components: (1) a game-theoretic MPC problem formulation that considers the interaction between CAVs and HDVs as a potential game, (2) an inverse learning method to learn online human driving model, and (3) a heuristic strategy to adapt the weights of the MPC problem.
We have also developed automatic adaptation strategies \cite{Le2023ACC,le2024controller} for the MPC problem based on Bayesian optimization and learning the solutions of contextual Bayesian optimization, respectively, which perform better than MPC with a heuristic strategy based on social value orientation in \cite{Le2022CDC} or non-adaptive MPC.

On the other hand, to the best of our knowledge, the coordination of multiple CAVs among HDVs remains a significant and underexplored challenge in the current research literature.
Formulating optimization problems for the coordination of multiple CAVs is much more challenging than for the control of single CAVs due to the following reasons.
First, the interaction between vehicles in such scenarios becomes more complex and difficult to model.
In addition, the control objective involves multiple traffic metrics, e.g., travel time and traffic smoothness, and some social metrics, such as human comfort and safety.
More importantly, the large number of CAVs and HDVs involved in the traffic scenarios leads to large-scale learning and optimization problems which might hinder the real-time applicability of the framework.
Recently, we presented a control framework that aims to derive time-optimal trajectories for CAVs in a mixed-traffic merging scenario given the HDVs’ future trajectories predicted from Newell’s car-following model \cite{Le2023CDC}.
The time-optimal trajectories are then combined with a safety filter based on control barrier functions.
We also proposed a stochastic time-optimal trajectory planning framework for coordinating multiple CAVs in mixed-traffic merging scenarios \cite{Le2023Stochastic}. 
To efficiently learn the driving behavior of human drivers online, we presented a data-driven model that combines Newell's car-following model with Bayesian linear regression. 
Using this prediction model and uncertainty quantification, we formulated a stochastic time-optimal control problem to find robust trajectories for CAVs. 
Additionally, we integrated a replanning mechanism to determine when deriving new trajectories for CAVs is necessary based on the accuracy of the Bayesian linear regression predictions. 
However, current methods analyze, design, and optimize a mobility system without considering the social dimension resulting in systems that might not be acceptable by the drivers, travelers, and the public. In the IDS Lab, we combine the three aforementioned dimensions \cite{Chremos2020MechanismDesign,chremos2020MobilityMarket,chremos2020SharedMobility,chremos2021MobilityGame}.

Our research has also expanded to autonomous mobility-on-demand (AMoD) systems. We have developed algorithms for routing with electric vehicles and charging scheduling \cite{bang2021AEMoD}. 
While both CAVs and AMoD systems fall under the umbrella of emerging mobility systems, their distinct objectives and perspectives became apparent.
Control of CAVs involves considering all physical phenomena, whereas AMoD systems revolve around optimization at a macroscopic level.
This disparity inspired my dissertation focus on bridging the gap between these realms and developing innovative solutions for the effective coordination and management of both CAVs and AMoD systems in emerging mobility scenarios.
To bridge the gap, we have developed a combined framework that considers the behavior and influence of CAVs in AMoD systems, optimizing routes for each CAV while addressing interactions of CAVs at intersections \cite{Bang2022combined,Bang2022rerouting}. 
This framework significantly enhanced traffic flow by predicting actual traffic conditions at intersections through a coordination problem. To ensure scalability, we extended this research by designing a hierarchical framework. This framework optimizes vehicle flow for routing at a macro level, extracts routes, assigns them to CAVs, and efficiently coordinates and controls CAVs at a micro level \cite{Bang2023flowbased}.
The approach offers superior computational efficiency by utilizing algorithmic methods, guaranteeing the achievability of optimal flow based on actual CAV movements. Additionally, we conducted a comprehensive study on the trade-offs associated with employing algorithmic approaches, determining the necessary conditions for optimal performance \cite{bang2023optimal}.
We have also employed micro-simulation to analyze CAVs' energy impact \cite{bang2023exploring} and modeled autonomous vehicles' strategic interactions using game theory \cite{chremos2023AtomicRouting}.
Recognizing the impact of emerging mobility systems on social inequities, we developed novel equity metrics to optimize traffic flow and achieve better equity \cite{Bang2023mem,bang2024cts, bang2024emergingequity}. These endeavors collectively contribute to the overarching goal of my research — developing comprehensive solutions for the efficient operation of emerging mobility systems.

A significant focus of our ongoing work is considering \textbf{human involvement in CPSs}. 
When humans are beneficiaries of autonomous decision-making, e.g., socio-technical systems, equity is a principal concern.
Recognizing the impact of emerging mobility systems on social inequities, we have developed a novel \textbf{mobility equity metric} for transportation networks \cite{bang2024cts}. %that uses only publicly available data to quantify the notion of equity across multiple modes of transportation.
%addresses the pressing need for \textbf{mobility equity metrics} in transportation networks.
%Equity must necessarily be a fundamental consideration in any CPS where humans are either participants or beneficiaries of autonomous decision-making.
%Recognizing the potential impact of emerging technologies on societal inequities, we've developed a novel metric that quantifies equity in multi-modal emerging mobility systems using only publicly available data.
We have utilized this metric to evaluate mobility equity using real data across U.S. cities \cite{bang2024emergingequity} and \textbf{optimize the flow of traffic} in emerging mobility systems \cite{Bang2023mem}.
In parallel, we have considered the explicit involvement of humans in the decision-making loop of Cyber-Physical Human Systems (CPHSs). 
In the context of autonomous driving, we have leveraged the framework reported in \cite{Nishanth2023AISmerging} and conformal prediction \cite{bang2024confidence}, to learn models for \textbf{human-driving} and generate\textbf{safe and effective control strategies} for automated vehicles merging among human driven vehicles.
For more general CPHSs, we have considered scenarios where an AI recommends decisions to a human, who may ignore recommendations to implement an action instinctively, leading to sub-optimal performance.
To address these challenges, we have proposed an \textbf{adherence-aware Q-learning} algorithm that outperforms baselines by adapting to the human's behavior \cite{faros2023adherence}. Building upon this result, we have also proposed the framework of human-AI POMDPs that allow us to learn more complex human behaviors such as dynamically evolving trust, and improve AI recommendations \cite{dave2024airecommend}.

}

\section{\textbf{Sequential Team Decision Proplems with Nonclassical Information Structures}}
\justify{

Team theory is a mathematical formalism for decentralized stochastic control problems  in which a ``team," consisting of a number of members, cooperates to achieve a common objective. It was developed to provide a rigorous mathematical framework of cooperating members in which all  members have the same objective yet different information.  The underlying structure to model a team decision problem consists of (1) a number of $K \in \mathbb{N}$ members of the team; (2) the decisions of each member; (3) the information available to each member, which is different; (4) an objective, which is the same for all members; and (5) the existence, or not, of communication between team members. Team theory can be applied effectively in applications that include informationally decentralized systems such as emerging mobility systems  \cite{zhao2019enhanced}, and in particular, optimal coordination of connected and automated vehicles at traffic scenarios \cite{Malikopoulos2017,Cassandras2017, mahbub2020decentralized,Malikopoulos2020, chalaki2020TCST}, networked control systems \cite{Hespanha:2007aa,Zhang:2020aa}, mobility markets \cite{chremos2020MobilityMarket}, smart power grids \cite{Khaitan:2013aa,Howlader:2014aa}, power systems \cite{Du:2017aa}, cooperative cyber-physical networks \cite{Pasqualetti:2015aa,Sami:2016aa,Clark:2017aa}, social media platforms \cite{Dave2020SocialMedia}, cooperation of robots \cite{Jadbabaie:2003aa,Saulnier:2017aa,Beaver2020AnFlockingb}, and internet of things \cite{Li:2016aa,Xu:2018aa,Ansere:2020aa}.	

Team theory was established with the seminal work of Marschak \cite{Marshak:1974aa}, Radner \cite{Radner1962}, and Marschak and Radner \cite{Marschak_Radner1972} on \textit{static team} problems,  and with Witsenhausen \cite{Witsenhausen1971,Witsenhausen1973} on \textit{dynamic team} problems. In static team problems \cite{Krainak:1982aa,Krainak:1982ab}, the information received by the team members is not affected by the decisions of other team members \cite{Yuksel2013}, while, in dynamic team problems, the information of at least one team member is affected by the decisions of other members  in the team\cite{Yuksel2013}. If there is a prescribed order in which team members make decisions, then such a problem is called a \textit{sequential} team problem. If, however, the team members make decisions in an order that depends on the realization of the team's uncertainty and decisions of other members, then such a problem is called a \textit{non-sequential} team problem. Formulating a well-posed non-sequential team problem is more challenging as we need to ensure that the problem is causal and deadlock free \cite{Andersland:1991aa, Andersland1992,Andersland1994}.

We have addressed sequential dynamic team decision problems with nonclassical information structures \cite{Dave2020,Dave2020a,Dave2021a,Dave2021nestedaccess,Dave2021minimax}. In the most recent effort \cite{Malikopoulos2021}, we provided structural results and a classical dynamic programming decomposition of sequential dynamic team decision problems. We first addressed the problem from the point of view of a manager who seeks to derive the optimal strategy for a team in a centralized process. Then, we addressed the problem from the point of view of each team member, and showed that their optimal strategies are the same as the ones derived by the manager. Our key contributions are (1) the structural results for the team from the point of view of a manager that yield an information state which does not depend on the control strategy of the team, and (2) the structural results for each team member that yield an information state which does not depend on their control strategy. These results allowed us to formulate two dynamic programming decompositions: (a) one for the team where the manager's optimization problem is over the space of the team's decisions, and (b) one for each team member where the optimization problem is over the space of the decision of each member. Finally, we showed that the control strategy of each team member is the same as the one derived by the manager. Therefore, each team member can derive their strategy, which is optimal for the team, without the manager's intervention.

A potential direction for future research should explore the intersection of learning and control \cite{Malikopoulos2022a,Malikopoulos2024} for team decision problems with nonclassical information structures. For example, cyber-physical systems, e.g., emerging mobility systems \cite{zhao2019enhanced}, in most instances, represent systems of systems with informationally decentralized structure. In such systems, however, there is typically a large volume of data with a dynamic nature which is added to the system gradually and not altogether in advance. Therefore, neither traditional supervised (or unsupervised) learning nor typical model-based control approaches can effectively facilitate feasible solutions with performance guarantees. These challenges could be circumvented at the intersection of learning and control. 

The IDS Lab's research activities in this area identify \textbf{key structural properties} of optimal strategies, reducing the computational burden of deploying control or reinforcement learning solutions \cite{dave2019decentralized, Dave2020a, Dave2021a, Dave2021nestedaccess, Dave2021minimax}. We have also extended our research to consider agents who compete to maximize their private utilities in \textbf{socio-technical systems}, leading to \textbf{strategic games}. 
More specifically, we proposed such a model for the interactions between social media platforms and a democratic government that seeks to \textbf{minimize misinformation} and designed a mechanism to incentivize misinformation filtering \cite{Dave2020SocialMedia}. Similarly, we studied the interactions among \textbf{travelers in an emerging mobility system} as a strategic game with bounded rationality of humans and studied traffic outcomes \cite{chremos2023AtomicRouting}.
%The resulting equilibria characterized the impact of autonomous vehicles on traffic conditions after rebound effects 
\\ \vspace{-9pt}
\\
We have also developed a theoretical framework introducing information states and facilitating \textbf{system approximations} for computationally efficient robust control \cite{Dave2022approx}.
The framework provides theoretical guarantees on the worst-case performance of approximate strategies, evaluated using only output data \cite{dave2022additive}, and thus, it enables both \textbf{supervised learning} and \textbf{reinforcement learning} for safety-critical CPSs \cite{Dave2023approximate, Dave2023infhorizon} such as automated vehicles and healthcare. This framework represents an important step towards a unification of learning and control methods for principled AI development.
We have also focused on utilizing theoretical foundations for CPS applications. In particular,  explored \textbf{CPSs vulnerable to attacks}, achieving an effective trade-off between optimal and robust strategies to avoid over-conservatism \cite{venkatesh2023stochastic}. %without succumbing to overly conservative strategies. In our preliminary results, we have extended my ideas to a constrained formulation that achieves a trade-off between optimal control and robust control. 
In simulation, we have used the approximate information state as an approximately sufficient statistic for predicting system behavior and decision-making by involving concepts of supervised learning and reinforcement learning. Subsequently, I adapted the notion of AIS to predict the future trajectory of a HDV for controlling a CAV in a mixed-traffic merging scenario \cite{Nishanth2023AISmerging}.

}	

\section{\textbf{Multi-Agent and Swarm System}}
\justify{

Flocking has the hallmark trait of emergence where agents are able to achieve a desired pattern at the system-level. 
In our recent review on optimal flocking \cite{Beaver2020AnFlockingb}, it became clear that a majority of flocking research focuses on implementations of Reynolds flocking rules.
In the IDS Lab, we are interested in how flocking emerges in swarms of constraint-driven agents, i.e., energy-minimizers subject to local interaction constraints  \cite{Beaver2020AnFlocking,Beaver2020BeyondFlocking,Beaver2021Constraint-DrivenStudy,bang2021energy}.
A critical aspect of this research is the development of self-relaxing constraints \cite{Beaver2020Energy-OptimalConstraints}, which can accommodate agents to re-plan their trajectories when new information renders their current state infeasible.
We have also proposed an approach that combines a heuristic banning mechanism with an energy-optimal local assignment \cite{Beaver2019AGeneration,Beaver2020AnAgents}, which hints at several interesting results. 

}

\newpage

\bibliographystyle{IEEEtran}
\bibliography{references1,IDS_Publications_05122024,TAC_Ref_structure,TAC_Ref_Andreas}

\end{document}